\documentclass[12pt]{amsart}

\usepackage{graphics}
\usepackage{amsmath}
\usepackage{amsfonts}
\usepackage{amssymb}
\usepackage{amscd}

\input xy
\xyoption{all}

\newtheorem{theorem}{\bf Theorem}

\newtheorem{corollary}[theorem]{\bf Corollary}

\title{Dessins d'enfants and differential equations}

\author{Finnur L\'arusson}
\address{School of Mathematical Sciences, University of Adelaide, Adelaide SA 5005, Australia.} 
\email{finnur.larusson@adelaide.edu.au}

\author{Timur Sadykov}
\address{Department of Mathematics and Computer Science, Krasnoyarsk State University, 660041
Krasnoyarsk, Russia.} 
\email{sadykov@lan.krasu.ru}

\thanks{Finnur L\'arusson was supported in part by the Natural Sciences and Engineering Research Council of Canada.  Timur Sadykov was supported in part by the Russian Foundation for Basic Research, grant 05-01-00517, by grant MK-851.2006.1 of the President of the Russian Federation, and by grant JSPS 06-01-91063 from the Japanese Society for the Promotion of Science.  Part of the work of both authors was done at the University of Western Ontario.}

\date{30 July 2006}

\begin{document}

\begin{abstract}
We state and solve a discrete version of the classical Riemann-Hilbert problem.  In particular, we associate a Riemann-Hilbert problem to every dessin d'enfants.  We show how to compute the solution for a dessin that is a tree.  This amounts to finding a Fuchsian differential equation satisfied by the local inverses of a Shabat polynomial.  We produce a universal annihilating operator for the inverses of a generic polynomial.  We classify those plane trees that have a representation by M\"obius transformations and those that have a linear representation of dimension at most two.  This yields an analogue for trees of Schwarz's classical list, that is, a list of those plane trees whose Riemann-Hilbert problem has a hypergeometric solution of order at most two.
\end{abstract}

\maketitle

\section{Review of the classical Riemann-Hilbert problem}

\noindent
Consider a homogeneous linear differential equation of order $n$ in one complex variable
$$y^{(n)}+p_1 y^{(n-1)} + \dots + p_{n-1}y' + p_n y = 0,$$
where $p_1,\dots,p_n$ are meromorphic at a point $a$ in the complex plane $\mathbb C$.  If $p_1,\dots,p_n$ are holomorphic at $a$, then $a$ is called a regular point of the equation; otherwise $a$ is called a singular point.  Cauchy proved that at a regular point, there is an $n$-dimensional space of holomorphic solution germs.

A singularity $a$ is called regular if every solution in a sector at $a$ grows at most polynomially at $a$.  This is equivalent to $a$ being Fuchsian, that is, to $p_i$ having a pole at $a$ of order at most $i$ for $i=1,\dots,n$.  (This is not so simple for the more general notion of first-order systems.)

Now we work globally on the Riemann sphere $\mathbb P$ and let $p_1,\dots,p_n$ be rational functions.  Let $S$ be the set of singularities.  Any holomorphic solution at a regular point $x_0$ can be analytically continued along any path in $\mathbb P\setminus S$.  Returning to $x_0$, we have another solution, linearly depending on the first one, and only depending on the homotopy class of the path.  This yields the monodromy representation
$$\pi_1(\mathbb P\setminus S, x_0) \to \text{Aut}(E) \cong \text{GL}(n,\mathbb C),$$
where $E$ is the $n$-dimensional solution space at $x_0$.

The equation is called Fuchsian if every point in $\mathbb P$ is regular or a regular singularity.  There is an extensive structure theory for regular singularities of first-order systems; irregular singularities are quite different.

The classical Riemann-Hilbert problem is as follows.  Given a finite set $S\subset\mathbb P$, a finite-dimensional complex vector space $E$, and a representation $\pi_1(\mathbb P\setminus S) \to \text{Aut}(E)$, is there a Fuchsian equation with this monodromy?  In other words, can we realize the elements of $E$ as holomorphic functions at a point in $\mathbb P \setminus S$, such that $E$ becomes the solution space of a homogeneous Fuchsian equation whose singularities lie in $S$ and the representation becomes the monodromy given by analytic continuation around the points of $S$?

By counting parameters, Poincar\'e saw that the answer is no.  We need to move to systems of first-order equations and realize the elements of $E$ as vector-valued functions.  In 1908, Plemelj gave an affirmative answer with a regular system.  For 70 years it was thought that his proof even gave a Fuchsian system, but in 1989, Bolibruch found examples where this is impossible.  Plemelj's argument works, however, when one of the monodromy matrices is diagonalizable.  There are several other positive results.  By the modern approach, begun by R\"ohrl in 1957, we can always solve the Riemann-Hilbert problem using a twisted Fuchsian system.  There is an extensive literature on the subject; we conclude this brief review by mentioning the survey \cite{Beauville} and the further references \cite{AnosovBolibruch}, \cite{GantzSteer}, and \cite{vanderPutSinger}.

\section{A discrete Riemann-Hilbert problem}

\noindent
Let $S\subset\mathbb P$ be a finite set with $m+1$ elements, $E$ be a finite set, and $F_m\cong\pi_1(\mathbb P \setminus S) \to \operatorname{Aut}(E)$ be a group homomorphism into the group of permutations of $E$.  Here $F_m$ denotes the free group on $m$ generators.  Assume the image acts transitively on $E$.  Can we realize the elements of $E$ as distinct holomorphic function germs at a point in $\mathbb P\setminus S$, such that $E$ spans the solution space of a homogeneous Fuchsian equation whose singularities lie in $S$ and the group homomorphism is the monodromy given by analytic continuation around the points of $S$?  We call this question the (discrete) Riemann-Hilbert problem associated to the given data.

By Plemelj's Theorem, we can realize the elements of $E$ by vector-valued holomorphic functions and solve the problem with a first-order Fuchsian system: we get a permutation representation on a vector space with basis $E$ and permutation matrices are diagonalizable.  There is a more elementary way of solving the problem with a single scalar equation that even represents elements of $E$ by algebraic functions.  The proof of the following result uses the well-known Wronskian construction of a differential equation with a prescribed solution space.

\begin{theorem}
Let $S$ be a finite subset of $\mathbb P$, $E$ be a finite set, and $\pi_1(\mathbb P \setminus S) \to \operatorname{Aut}(E)$ be a group homomorphism with a transitive image.   Then the elements of $E$ correspond to distinct algebraic functions that span the solution space of a homogeneous Fuchsian equation with singularities in $S$, such that the group homomorphism corresponds to the monodromy given by analytic continuation around the points of $S$.
\end{theorem}

\begin{proof}  
By basic covering space theory, the data of the discrete Riemann-Hilbert problem correspond to an unbranched covering space over $\mathbb P\setminus S$ with generic fibre $E$ and the image of $\pi_1(\mathbb P \setminus S)$ as its monodromy group (not to be confused with the group of covering transformations).  Such a covering space has a unique complex structure making the covering map holomorphic and it extends uniquely to an $n$-sheeted branched covering space $p:X\to\mathbb P$, where $n=|E|$.  Transitivity means that $X$ is connected.

Let $\sigma_1,\dots,\sigma_n$ be the germs of the inverses of $p$ at a chosen base point $x_0$ in $\mathbb P\setminus S$.  Find a meromorphic function $h$ on the compact Riemann surface $X$ separating the distinct points $\sigma_1(x_0),\dots,\sigma_n(x_0)$.  Then the germs $s_i=h\circ\sigma_i$, $i=1,\dots,n$, are distinct and there is a natural bijection
$E\to\{s_1,\dots,s_n\}$ such that the given permutations of $E$ correspond to analytic continuation around the points of $S$.  Also, $s_1,\dots,s_n$ are algebraic functions.  In fact, 
$$s_i^n+c_1 s_i^{n-1}+\dots+c_{n-1} s_i+c_n=0,$$
where the rational functions $c_1,\dots,c_n$ on $\mathbb P$ are the elementary symmetric functions of the values of $h$ on the fibres of $p$.

Now $s_1,\dots,s_n$ need not be linearly independent; suppose $s_1,\dots, s_k$ form a basis for the vector space spanned by them.  Then this vector space is the solution space of the $k^{\text{th}}$-order homogeneous monic linear differential equation
$$(-1)^{k+1}\left| \begin{matrix}
y & s_1 & \hdots & s_k \\
y' & s_1' & \hdots & s_k' \\
\vdots & \vdots &  & \vdots \\
y^{(k)} & s_1^{(k)} & \hdots & s_k^{(k)}
\end{matrix} \right|
\left| \begin{matrix}
s_1 & \hdots & s_k \\
s_1' & \hdots & s_k' \\
\vdots &  & \vdots \\
s_1^{(k-1)} & \hdots & s_k^{(k-1)}
\end{matrix} \right|^{-1} = 0.
$$
It is easy to check that the coefficients of this equation are invariant under monodromy, so they extend holomorphically to $\mathbb P\setminus S$.  Also, since $s_1,\dots,s_k$ are algebraic functions, so are the coefficients, so they are rational functions.  Moreover, since any solution is a linear combination of the algebraic functions $s_1,\dots,s_k$, any solution in a sector at a point of $S$ grows at most polynomially, so the equation is regular and hence Fuchsian.
\end{proof}

\section{The case of three singularities: dessins d'enfants}

\noindent
From now on, we restrict ourselves to the case when $S$ has (at most) three points, say $S=\{0,1,\infty\}$.  This case is of particular interest for two reasons.  First, the data of a finite set $E$ and a group homomorphism $F_2 \to \operatorname{Aut}(E)$ determines the finite branched covering $p:X\to\mathbb P$ up to isomorphism, whereas for $|S|=m+1>3$, the covering depends on $m-2$ continuous parameters, corresponding to the positioning of the points of $S$ modulo M\"obius transformations.  Second, the data can be encoded combinatorially as a dessin d'enfants, a graph with a certain extra structure, providing a new and fruitful point of view.

Furthermore, there is a profound link with number theory by a theorem of Belyi \cite{Belyi1, Belyi2}, stating that the finite branched covering spaces of $\mathbb P$, branched over at most three points, are precisely the compact Riemann surfaces that can be defined over the field of algebraic numbers, that is, over some number field.  The absolute Galois group $\operatorname{Aut}(\overline{\mathbb Q}/\mathbb Q)$ thus acts, in fact faithfully, on the set of dessins.  Grothendieck had the idea, expounded in his {\it Esquisse d'un programme} (1984, \cite{Grothendieck}), to use this action to study the absolute Galois group.  The main question is whether the Galois orbits of dessins can be characterized by combinatorial or topological invariants.  Despite considerable work, this question is still very much open.

The data of the discrete Riemann-Hilbert problem is two permutations acting transitively on a finite set $E$.  Call one permutation black and the other white.  Construct a graph by taking the orbits (cycles) of each of the two permutations as vertices and the elements of $E$ as edges.  Transitivity means that the graph is connected.  Each edge joins a black vertex and a white vertex, so the graph is bicoloured.  Finally, there is a cyclic ordering of the edges at each vertex.  This is a dessin d'enfants.  The reader who is new to dessins might consult the survey \cite{Oesterle} or the books \cite{LandoZvonkin} or \cite{Schneps}.

Dessins d'enfants correspond bijectively to the data of two permutations acting transitively on a finite set, and to finite branched covering spaces of $\mathbb P$, branched over 0, 1, and $\infty$ (up to isomorphisms).  Such covering maps are called Belyi maps.  The dessin sits in the covering space as the preimage of the segment $[0,1]$.  Its complement is a union of cells, or faces, each containing one preimage of $\infty$.  If the dessin has $v$ vertices, $e$ edges, and $f$ faces, then $v-e+f$ equals the Euler characteristic of the covering space.

So we have posed and, in Theorem 1, solved a Riemann-Hilbert problem for dessins.  Various questions arise.  Can the problem be solved explicitly or algorithmically, say for special classes of dessins such as trees?  For which dessins is there a solution by a Fuchsian equation of a special type, say a second-order equation or a hypergeometric equation?  These questions will be addressed in what follows.  More ambitiously, we could ask whether relating Fuchsian equations to dessins might lead to new and useful dessin invariants.  This is completely open.

\section{Explicitly solving the Riemann-Hilbert problem for a tree}

\noindent
We shall discuss the problem of actually computing the solution of the Riemann-Hilbert problem for a dessin.  Clearly, calculating the inverses that go into the Wronskians in the solution given above is very hard in general, so another way must be found.  We will restrict ourselves to the case of a tree (a connected graph without circuits), that is, to the case when $v=e+1$.  Then Euler's formula implies that, first, $g=0$, so the covering space is $\mathbb P$.  (Calculations for dessins of higher genus turn out to be quite hard.)  Second, $f=1$, so the preimage of $\infty$ is a single point, which we can take to be $\infty$ itself, so the Belyi map is a polynomial with at most two critical values.  Such a polynomial is called a Shabat polynomial.  

The correspondence between dessins and Belyi maps specializes as follows.  Shabat polynomials $P$ and $Q$ are said to be equivalent if there are invertible affine transformations (holomorphic automorphisms of $\mathbb C$) $\phi$ and $\psi$ such that $\phi\circ P=Q\circ\psi$.  A plane tree is a tree with a cyclic ordering of the edges at each vertex.  There is a bijective correspondence between equivalence classes of Shabat polynomials and isomorphism classes of plane trees.  (Here, we do not distinguish the two critical points of a Shabat polynomial or the two bicolourings of a tree.)  For an introduction to these notions, see \cite{BetremaZvonkin}.

\smallskip \noindent
{\it Simplest examples.}  (1)  Shabat polynomials with only one critical value, just $x^n$ modulo affine coordinate changes, for which the tree is a star with $n$ arms, depicted here for $n=4$.  
$$
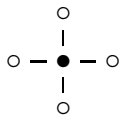
\xymatrix@=2mm{ & {\circ} \ar@{-}[d] & \\ {\circ} \ar@{-}[r] & {\bullet} \ar@{-}[d] \ar@{-}[r] & {\circ} \\ & {\circ}}$$
The inverses $x^{1/n}$ satisfy the first-order (hypergeometric) equation $y'-\frac 1{nx}y=0$.

(2)  The Chebyshev polynomials $\cos(n\arccos x)$, for which the tree is a chain with $n$ edges.  By a chain we mean a tree with two ends; then all the other vertices have valency two.  The second-order hypergeometric equation 
$$(1-x^2)y''-xy'+\tfrac 1{n^2}y=0$$
for the inverses $\cos(\tfrac 1 n\arccos x)$ is well known.  The affine change of variables $t=(1+x)/2$ transforms this equation into
$$t(1-t)y''+(\tfrac 1 2 - t)y'+\tfrac1{n^2}y=0,$$
which is a special case of Gauss' hypergeometric equation in its traditional form
$$x(1-x)y''+(c-x(a+b+1))y'-aby=0.$$

(3) The Shabat polynomials $(x^m-1)^2$, for which the tree is a star with $m$ arms, each composed of two edges, depicted here for $m=4$.  
$$
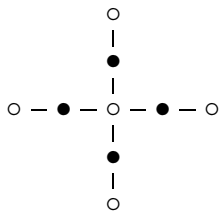
\xymatrix@=2mm{ & & {\circ} \ar@{-}[d] & & \\ & & {\bullet} \ar@{-}[d] & & \\ {\circ} \ar@{-}[r] & {\bullet} \ar@{-}[r] & {\circ} \ar@{-}[r] \ar@{-}[d] & {\bullet} \ar@{-}[r] & {\circ} \\ & & {\bullet} \ar@{-}[d] & & \\ & & {\circ}}$$
Let us call such a tree a 2-star.  The inverses $(x^{1/2}+1)^{1/m}$ satisfy the hypergeometric equation
$$x(1-x)y''+((\tfrac 1 m-\tfrac 3 2)x + \tfrac 1 2)y'+\tfrac 1{4m}(1-\tfrac 1 m)y=0.$$

\smallskip

It turns out that for a dessin which is a tree, both the Belyi map and a Fuchsian equation that solves the Riemann-Hilbert problem can be explicitly calculated.  The former is well known (see e.g.\ \cite{LandoZvonkin}, sec.\ 2.2).  Take a bicoloured tree in $\mathbb P$ with $p$ black vertices $a_1,\dots,a_p$ with valencies $\alpha_1,\dots,\alpha_p$, and $q$ white vertices $b_1,\dots,b_q$ with valencies $\beta_1,\dots,\beta_q$.  Then $\sum\alpha_i=\sum\beta_j=e$.  The tree property is that $p+q=e+1$.  The corresponding Shabat polynomial $P$, normalized to have critical values in $\{0,1\}$, has degree $n=e$ and satisfies the equations
$$\begin{aligned}
P(x)&=c(x-a_1)^{\alpha_1}\dots(x-a_p)^{\alpha_p}, \\
P(x)-1&=c(x-b_1)^{\beta_1}\dots(x-b_q)^{\beta_q}.
\end{aligned}$$
Given $p$, $q$, and the valencies of an abstract tree, this set of equations has only finitely many solutions $a_1,\dots,a_p,b_1,\dots,b_q$ up to normalization, and these solutions are computable in practice.  One of them, unique up to equivalence, is the corresponding Shabat polynomial.  The Shabat property is that $p+q$ equals $n+1$, which is the smallest $p+q$ can be for a polynomial $P$.  Equivalently,
$$P'(x)=cn(x-a_1)^{\alpha_1-1}\dots(x-a_p)^{\alpha_p-1}(x-b_1)^{\beta_1-1}
\dots(x-b_q)^{\beta_q-1}.$$

The solution of the Riemann-Hilbert problem for the tree is a Fuchsian equation of order at most $n$ with polynomial coefficients $q_0,\dots,q_n$ such that $q_n\sigma^{(n)}+\dots+q_0\sigma=0$ for any or all local inverses $\sigma$ of $P$ at a regular value.  Note that if $n\geq 2$ and the penultimate coefficient of $P$ is nonzero, then $q_0=0$; otherwise, the inverses span a space of dimension less than $n$ and we can take $q_n=0$.

Precomposing the equation by $P$ gives $\sum (q_k\circ P)r_k=0$, where $r_k=\sigma^{(k)}\circ P$.  It is straightforward that $r_k$ is a rational function, in fact a polynomial in $P',\dots,P^{(k)}$ divided by $(P')^{2k-1}$, and it can be easily calculated using the recursive formulas
$$r_0=x, \qquad r_{k+1}=\sigma^{(k+1)}\circ P = \frac{(\sigma^{(k)}\circ P)'}{P'}=\frac{r_k'}{P'}.$$
Also, $r_k$ can be expressed by an integral formula and as a residue.  Thus, with some work, $q_0,\dots,q_n$ can be computed by converting the equation to a linear system.

The simplest nontrivial example is the following tree, the smallest tree that is neither a star nor a chain.  
$$
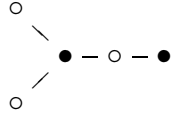
\xymatrix@=2mm{{\circ} \ar@{-}[dr] & & &  \\ & {\bullet} \ar@{-}[r] \ar@{-}[dl] & {\circ} \ar@{-}[r] & {\bullet} \\ {\circ}}$$
For this tree, we can take $P(x)=4x^3-x^4$, and the associated Fuchsian equation
$$(32x^4-864x^3)y'''' + (208x^3-3456x^2)y''' + (270x^2-1920x)y'' + 45xy' =  0$$
is in fact hypergeometric. 

\section{The annihilator of the inverse of a generic polynomial}

\noindent
Let $P$ be a complex polynomial of degree $n\geq 1$, viewed as a finite branched covering map $\mathbb P\to\mathbb P$.  Let $\sigma_1,\dots,\sigma_n$ be the germs of the inverses of $P$ at a chosen regular value $x_0$ of $P$.  Say $\sigma_k(x_0)=t_k$.  Taking a small parametrized circle $\gamma$ around $t_k$, it is well known that for $x$ inside $P\circ\gamma$, 
$$\sigma_k(x)=\frac 1{2\pi i}\int_\gamma\frac{tP'(t)}{P(t)-x}\,dt.$$
Differentiating $m\geq 1$ times yields
$$\begin{aligned}\sigma_k^{(m)}(x)&=\frac {m!}{2\pi i}\int_\gamma\frac{tP'(t)}{(P(t)-x)^{m+1}}\,dt
=(m-1)!\operatorname{Res}_{\sigma_k(x)}\frac{mtP'(t)}{(P(t)-x)^{m+1}} \\
&=(m-1)!\operatorname{Res}_{\sigma_k(x)}(P(t)-x)^{-m}, \end{aligned}$$
since
$$\frac d{dt}\frac t{(P(t)-x)^m} = \frac 1{(P(t)-x)^m} - \frac{mtP'(t)}{(P(t)-x)^{m+1}}$$
has no residue at $\sigma_k(x)$.  By a familiar formula for residues at higher-order poles (see \cite{Remmert}, p.~273),
$$\sigma_k^{(m)}(x)=\frac{d^{m-1}}{dt^{m-1}}\bigg(\frac{t-\sigma_k(x)}{P(t)-x}\bigg)^m\bigg|_{t=\sigma_k(x)}.$$
Here, the fraction is to be holomorphically continued across $\sigma_k(x)$.  Now, assuming $P$ is monic,
$$P(t)-x=(t-\sigma_1(x))\cdots(t-\sigma_n(x)),$$
so
$$\frac{P(t)-x}{t-\sigma_k(x)}= \prod_{i\neq k}(t-\sigma_i(x)).$$
We conclude that in a neighbourhood of $x_0$, and thus wherever $\sigma_1,\dots,\sigma_n$ are all defined, we have
$$\sigma_k^{(m)}=Q_{m,k}(\sigma_1,\dots,\sigma_n),$$
where $Q_{m,k}$ is the rational function
$$Q_{m,k}(t_1,\dots,t_n)=\frac{d^{m-1}}{dt^{m-1}}\prod_{i\neq k}(t-t_i)^{-m}\bigg|_{t=t_k}, \qquad 1\leq k\leq n,\quad m\geq 1.$$
For $k=0,\dots,n$, the coefficient of $y^{(k)}$ in the determinant
$$ \left| \begin{matrix}
y & t_1 & \hdots & t_n \\
y' & Q_{11}(t_1,\dots,t_n) & \hdots & Q_{1n}(t_1,\dots,t_n) \\
\vdots & \vdots & & \vdots \\
y^{(n)} & Q_{n1}(t_1,\dots,t_n) & \hdots & Q_{nn}(t_1,\dots,t_n)
\end{matrix} \right| $$
is a symmetric rational function in $t_1,\dots,t_n$, so it is a rational function $R_k$ in the elementary symmetric functions in $t_1,\dots,t_n$.  Writing 
$$P(t)=t^n + a_{n-1} t^{n-1} + \dots + a_1 t + a_0,$$
the elementary symmetric functions in $\sigma_1(x),\dots,\sigma_n(x)$ are
$-a_{n-1}=\sigma_1(x)+\dots+\sigma_n(x)$, $a_{n-2}, \dots, (-1)^{n-1}a_1$, $(-1)^{n}(a_0-x)=\sigma_1(x)\cdots\sigma_n(x)$.  Hence, assuming the inverses $\sigma_1,\dots,\sigma_n$ are linearly independent, they satisfy the differential equation
$$\sum_{k=0}^n R_k\big(-a_{n-1},\dots,(-1)^{n-1}a_1, (-1)^n(a_0-x)\big)\frac{d^k y}{dx^k} =0.$$
Thus we have a universal annihilating operator with polynomial coefficients for the inverses of a generic monic polynomial, computable, although with rapidly increasing labour as the degree grows.  It is the same operator as in the proof of Theorem 1, except that there it was made monic by dividing by the leading coefficient.  The crucial difference is that the construction here does not require us to calculate the inverses $\sigma_1,\dots,\sigma_n$ themselves.  Note that since the inverses are linearly independent by assumption, if $n\geq 2$, their sum is a nonzero constant solution, so the term of order zero vanishes.

\smallskip \noindent
{\it Examples.}  (1)  The roots of the monic quadratic equation $y^2+a_1y+a_0=0$ are annihilated by the differential operator
$$D_2 = (a_1^2 - 4a_0)\frac{\partial^2}{\partial a_0^2} - 2 \frac \partial {\partial a_0}.$$
To obtain the annihilating operator for the inverses of $t^2+a_1 t+a_0$, simply replace $a_0$ by $a_0-x$.  This comment also applies to the following examples.

(2)  The roots of the cubic equation $y^3+a_2 y^2+a_1y+a_0=0$ are annihilated by the differential operator
$$D_3 =(a_1^2 a_2^2 - 4a_0 a_2^3 - 4a_1^3 + 18 a_0a_1a_2 - 27a_0^2) \frac{\partial^3}{\partial a_{0}^3} -
(6a_2^3 - 27 a_1 a_2 + 81a_0)\frac{\partial^2}{\partial a_0^2} - 24\frac{\partial}{\partial a_0}.$$

(3)  Finally, the roots of the quartic equation $y^4+a_3y^3+a_2 y^2+a_1y+a_0=0$ are annihilated by the differential operator
\begin{multline*}D_4 = (45 a_1^2 + 8 a_0 a_2 + 14 a_2^3  - 47 a_1 a_2 a_3- 3 a_0 a_3^2  - 4 a_2^2 a_3^2 + 12 a_1 a_3^3) \\ (-256 a_0^3 + 27 a_1^4 - 144 a_0 a_1^2 a_2 +128 a_0^2 a_2^2 + 4 a_1^2 a_2^3 - 16 a_0 a_2^4 +192 a_0^2 a_1 a_3 - 18 a_1^3 a_2 a_3 \\ + 80 a_0 a_1 a_2^2 a_3 + 6 a_0 a_1^2 a_3^2 - 144 a_0^2 a_2 a_3^2 - a_1^2 a_2^2 a_3^2+ 4 a_0 a_2^3 a_3^2 + 4 a_1^3 a_3^3 - 18 a_0 a_1 a_2 a_3^3 + 27 a_0^2 a_3^4) \frac{\partial^4}{\partial a_0^4} \\ + 4 (-21600 a_0^2 a_1^2 - 3328 a_0^3 a_2 - 4104 a_1^4 a_2 +6768 a_0 a_1^2 a_2^2 - 5696 a_0^2 a_2^3 - 1718 a_1^2 a_2^4 + 2192 a_0 a_2^5 \\ - 140 a_2^7 + 10800 a_0 a_1^3 a_3 + 24096 a_0^2 a_1 a_2 a_3 +6516 a_1^3 a_2^2 a_3 - 3920 a_0 a_1 a_2^3 a_3 +1170 a_1 a_2^5 a_3 + 1248 a_0^3 a_3^2 \\ + 189 a_1^4 a_3^2 -19200 a_0 a_1^2 a_2 a_3^2 + 384 a_0^2 a_2^2 a_3^2 -1820 a_1^2 a_2^3 a_3^2 -3130 a_0 a_2^4 a_3^2 + 75 a_2^6 a_3^2 -6336 a_0^2 a_1 a_3^3 \\ - 1784 a_1^3 a_2 a_3^3 +9276 a_0 a_1 a_2^2 a_3^3 - 595 a_1 a_2^4 a_3^3 +4392 a_0 a_1^2 a_3^4 + 648 a_0^2 a_2 a_3^4 +1113 a_1^2 a_2^2 a_3^4 + 1188 a_0 a_2^3 a_3^4 \\ -10 a_2^5 a_3^4 + 48 a_1^3 a_3^5 -3726 a_0 a_1 a_2 a_3^5 + 75 a_1 a_2^3 a_3^5 -81 a_0^2 a_3^6 - 135 a_1^2 a_2 a_3^6 -135 a_0 a_2^2 a_3^6 + 405 a_0 a_1 a_3^7)\frac{\partial^3}{\partial a_0^3} \\ +60 (-2244 a_0 a_1^2 - 288 a_0^2 a_2 + 389 a_1^2 a_2^2 -656 a_0 a_2^3 + 118 a_2^5 + 561 a_1^3 a_3 +2388 a_0 a_1 a_2 a_3\\ -a_1 a_2^3 a_3 + 108 a_0^2 a_3^2 -1011 a_1^2 a_2 a_3^2 + 141 a_0 a_2^2 a_3^2 -165 a_2^4 a_3^2 - 615 a_0 a_1 a_3^3 + a_1 a_2^2 a_3^3\\ +228 a_1^2 a_3^4 +24 a_0 a_2 a_3^4 + 62 a_2^3 a_3^4 -195 a_1 a_2 a_3^5 - 3 a_0 a_3^6 - 7 a_2^2 a_3^6 +21 a_1 a_3^7)\frac{\partial^2}{\partial a_0^2} \\ -120 (243 a_1^2 + 24 a_0 a_2 + 74 a_2^3 -249 a_1 a_2 a_3 - 9 a_0 a_3^2 - 21 a_2^2 a_3^2 +63 a_1 a_3^3)\frac{\partial}{\partial a_0}.\end{multline*}
From this formula we can recover the solution of the Riemann-Hilbert problem for the tree given at the end of the previous section.  Whereas the leading coefficients of $D_2$ and $D_3$ are the discriminants of the corresponding equations, the leading coefficient of $D_4$ is the product of the discriminant and the polynomial $45a_1^2+\dots+12a_1 a_3^3$.

\section{M\"obius trees}

\noindent
Dessins can be represented in all sorts of structures.  Any time we have a set $S$ with additional structure, we can consider the possibility of realizing the edges of a dessin as distinct points of $S$ and the permutations of the dessin as automorphisms of the structure.  In this section we shall describe those plane trees that can be represented in the Riemann sphere $\mathbb P$ with its group $\operatorname{Aut}\mathbb P$ of M\"obius transformations.  We call such trees M\"obius trees.  This result will be used in the following section to describe those trees whose Riemann-Hilbert problem has a solution of order at most two.

It is instructive, although not necessary for the proof of our theorem, to keep in mind the classification of finite M\"obius groups (see \cite{Maskit}, Sec.~V.C).  Every finite subgroup of $\operatorname{Aut}\mathbb P$ is conjugate to a cyclic group generated by a rotation $z\mapsto \theta z$, where $\theta$ is a primitive $n$-th root of unity; a dihedral group of order $2n$ generated by such a rotation and the inversion $z\mapsto 1/z$; or one of three sporadic groups, the orientation-preserving symmetry groups of the Platonic solids.  The symmetry group of the tetrahedron has order 12, that of the cube and the octahedron order 24, and that of the icosahedron and the dodecahedron order 60.  

If $T$ is a M\"obius tree with at least three edges, then its two permutations extend uniquely to M\"obius transformations that generate a finite M\"obius group $G$.  Using the fact that a tree must have an end, that is, a fixed point for one of the permutations, it is easy to see that if $G$ is cyclic, then $T$ is a star.  If $G$ is dihedral and both permutations are half-turns $z\mapsto \theta^j/z$, then, since half-turns have order two, every non-end vertex has valency two, and $T$ is a chain.  If $G$ is dihedral and one permutation is a rotation and the other is a half-turn, then either an edge is fixed by the rotation, so the tree has only two edges, or there are only two ends, so the tree is a chain.  

Thus, a M\"obius tree associated to a cyclic or a dihedral group is either a star or a chain.  The following theorem shows that the sporadic groups contribute nothing further.  As remarked above, the proof does not use the classification of finite M\"obius groups.

\begin{theorem}  The M\"obius trees are precisely the stars and the chains.
\end{theorem}

In other words, a M\"obius tree either has as few ends as a tree can possibly have (two), or as many ends as a tree can possibly have (every vertex except one is an end).

\begin{proof}  Clearly, a star has a M\"obius representation with one permutation being a rotation and the other the identity.  Also, a chain with $n$ edges has a M\"obius representation with permutations $z\mapsto 1/z$ and $z\mapsto \theta/z$, where $\theta$ is a primitive $n$-th root of unity, and the set of edges is the orbit of $1$.

Conversely, say $T$ is a M\"obius tree with $e$ edges and $v=e+1$ vertices.  Call its permutations $A$ and $B$.  A M\"obius transformation of finite order is conjugate to a rotation, so if it is not the identity, it has exactly two fixed points and the orbits of all other points have the same size.  Say $A$ has $\epsilon$ ends and $\alpha$ cycles of size $a\geq 2$, and $B$ has $\eta$ ends and $\beta$ cycles of size $b\geq 2$.  Then 
$$e=\epsilon+a\alpha = \eta+b\beta, \qquad e+1=v=\epsilon+\alpha+\eta+\beta.$$
Assume $T$ is not a star, so $\alpha,\beta\geq 1$ and neither $A$ nor $B$ is the identity, so $\epsilon,\eta\leq 2$.  Also, assume $T$ is not a chain, so its number of ends is $\epsilon+\eta\geq 3$; furthermore, $a\geq 3$ or $b\geq 3$.  Thus $e\geq 4$, so the action of $A$ and $B$ on the set of edges determines them as M\"obius tranformations.  They generate a finite M\"obius group $G$.

Since $T$ has at least three ends and $A$ and $B$ have only two fixed points each, both $A$ and $B$ fix at least one edge each.  Consider an edge fixed by $A$, viewed as a point $p\in\mathbb P$.  If the stabilizer of $p$ in $G$ is larger than the subgroup generated by $A$, then there is a word in $A$ and $B$, reduced in $G$, in which $B$ does occur, that fixes $p$, but this yields a circuit in $T$.  Thus the size of the stabilizer of $p$ is the order $a$ of $A$.  Similarly, the stabilizer of an edge fixed by $B$ has size $b$.  Since $G$ acts transitively on the edges, all stabilizers have the same size, so $a=b$.

Now $T$ has three or four ends.  We will eliminate both cases.  If $T$ has three ends, say $\epsilon=1$ and $\eta=2$ (the other case, $\epsilon=2$ and $\eta=1$, is analogous), then
$$\alpha+\beta+2=e=1+a\alpha=2+b\beta=2+a\beta$$
and $a(a-2)\beta=1$, which is impossible.

If $T$ has four ends, so $\epsilon=\eta=2$, then 
$$a\beta=b\beta=a\alpha=e-2=\alpha+\beta+1,$$
so $\alpha=\beta$ and $(a-2)\alpha=1$.  The only solution is $\alpha=\beta=1$, $a=b=3$, so $T$ is the following tree (or its opposite colouring).
$$
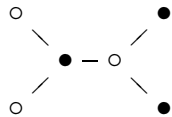
\xymatrix@=2mm{{\circ} \ar@{-}[dr] & & & {\bullet} \ar@{-}[dl]  \\ & {\bullet} \ar@{-}[r] \ar@{-}[dl] & {\circ} \ar@{-}[dr] \\ {\circ} & & & {\bullet}}$$
This tree is not a M\"obius tree.  Namely, after a suitable conjugation we have $A(z)=\theta z$, where $\theta=e^{2\pi i/3}$.  The five edges of the tree are represented by the fixed points $0$ and $\infty$ of $A$, as well as $c$, $\theta c$, and $\theta^2 c$ for some $c\in\mathbb C^\times$.  Say the three-element $B$-orbit contains $0$, $\infty$, and $c$, such that $B(0)=\infty$, $B(\infty)=c$, and $B(c)=0$.  These conditions determine $B$.  We get $B(z)=c(z-c)/z$, but then $B$ does not fix $\theta c$ or $\theta^2 c$.
\end{proof}

\section{Trees of order at most two}

\noindent
An $n$-dimensional linear representation of a dessin realizes its edges as distinct vectors in an $n$-dimensional (here, complex) vector space and the permutations of the dessin as linear automorphisms.  Every dessin has a linear representation by permutation matrices in a space of dimension equal to its number of edges.  We are concerned here with linear representations of small dimension and we shall restrict our attention to trees.  It is easy to see that the plane trees that have a one-dimensional linear representation are precisely the stars.  The purpose of our final section is to prove the following result.

\begin{theorem}  Let $T$ be a plane tree.  The following are equivalent.
\begin{enumerate}
\item[(i)]  $T$ has a linear representation of dimension at most two. 
\item[(ii)] The Riemann-Hilbert problem for $T$ has a solution of order at most two. 
\item[(iii)]  The Riemann-Hilbert problem for $T$ has a hypergeometric solution of order at most two.
\item[(iv)]  $T$ is a star, a 2-star, or a chain.
\end{enumerate}
\end{theorem}

Thus, the analogue for trees of Schwarz's classical list contains three infinite families and no sporadics.

\begin{proof}  Clearly, (iii) $\Rightarrow$ (ii) $\Rightarrow$ (i).  If $T$ is a star, a 2-star, or a chain, then the equivalence class of Shabat polynomials of $T$ contains a monomial, $(x^m-1)^2$, or a Chebyshev polynomial, respectively.  From Section 4, we know that the inverses of any such polynomial satisfy a homogeneous hypergeometric equation of order at most two.  Thus, (iv) implies (iii).

We need to show that (i) implies (iv).  Let $T$ be a plane tree with a linear representation of dimension two.  Assume $T$ is nontrivial, that is, has more than one edge, so its edges are represented by points in $\mathbb C^2\setminus\{0\}$.  Pass to the projectivized dessin $S$, that is, identify those edges of $T$ that lie in the same line through the origin in $\mathbb C^2$.  The two permutations of $T$ induce permutations of the edges of $S$ by M\"obius transformations.  We obtain an epimorphism $T\to S$ inducing a factorization $P=Q\circ\phi$ of any Shabat polynomial $P$ associated to $T$.  The source of $Q$ is $\mathbb P$ since the genus does not increase, and $Q$ has only one point above $\infty$ like $P$, so $Q$ is a Shabat polynomial and $S$ is a M\"obius tree.  The critical points of $P$ are critical points of the finite branched covering $\phi:\mathbb P\to\mathbb P$ or points taken by $\phi$ to critical points of $Q$.

Now we must analyze how $T$ sits above $S$.  The critical values of $\phi$ are those vertices of $S$ that have vertices of $T$ of higher valency above them.  If $\phi$ has no critical values, then $\phi$ is an automorphism and $T$ is a M\"obius tree.  If an orbit of one of the permutations shrinks under projectivization, it has at least two edges in a line through the origin in $\mathbb C^2$, so the line is invariant under the permutation.  Thus the orbit is contained in the line and collapses to a single edge.  Hence, critical vertices of $S$ are ends.  

If $S$ has two distinct critical ends, then above one of them we must have both an end and a non-end, for if we had only vertices of valency at least two over both of them, then by going back and forth between the preimages of the critical ends we could create a circuit in $T$.  In this case, a line through the origin contains a fixed point of one permutation and a nontrivial orbit of the other, so the line is invariant under both permutations, the representation is one-dimensional, and $T$ is a star.

We have arrived at the case when $S$ has precisely one critical end.  Then $T$ is obtained from $S$ by gluing together $\deg\phi=\deg P/\deg Q\geq 2$ copies of $S$ at the critical end.  If $S$ is a star with at least three arms, then one permutation of $T$ will always map a fixed point for the other to another such point.  If the two fixed points are linearly dependent, then the representation is one-dimensional; otherwise, the second permutation is the identity: in either case $T$ would be a star. 

If $k=\deg\phi\geq 3$ and $S$ is a chain with at least three edges, then one of the permutations, call it $A$, has edge orbits of length two and $k$.  If the eigenvalues of $A$ are primitive roots of unity of orders $a$ and $b$, then the possible lengths of orbits of $A$ are $a$, $b$ --- these orbits lie in eigenlines of $A$ --- and the least common multiple of $a$ and $b$.  Hence, an orbit of length two lies in an eigenline of $A$, and the other permutation will have a fixed vector $v$ (an end) in such an orbit.  The line spanned by $v$ is then invariant under both permutations, so the representation is one-dimensional and $T$ is a star, which is absurd.  Thus, $T$ is obtained by gluing together two copies of a chain of arbitrary length, or any number of copies of a chain of length one or two, so (iv) holds.
\end{proof}

Finally, we use Theorem 3 to describe those Shabat polynomials that, modulo translation, have hypergeometric inverses of order at most two.  These are precisely the Shabat polynomials that are equivalent to one of the examples in Section 4.

\begin{corollary}  Let $P$ be a Shabat polynomial.  There is $b\in\mathbb C$ such that the inverses of the equivalent Shabat polynomial $P(x+b)$ satisfy a homogeneous Fuchsian or, equivalently, hypergeometric equation of order at most two if and only if the tree of $P$ is a star, a 2-star, or a chain.
\end{corollary}

\begin{proof}  If there is a Shabat polynomial equivalent to $P$ whose inverses satisfy a homogeneous Fuchsian equation of order at most two, then the Riemann-Hilbert problem for the tree of $P$ has a solution of order at most two, so by Theorem 3, the tree of $P$ is a star, a 2-star, or a chain.

If the tree of $P$ is a star, a 2-star, or a chain, then $P$ is equivalent to a monomial,  $(x^m-1)^2$, or a Chebyshev polynomial, respectively, so there are affine transformations $\phi:x\mapsto ax+b$ and $\psi:x\mapsto cx+d$, $a,c\neq 0$, such that the inverses $f$ of the equivalent polynomial $\psi\circ P\circ\phi$ satisfy a hypergeometric equation $y''+ry'+sy=0$.  We have $f(x)=g(x/c-d/c)/a$, where $g$ is an inverse of $P(x+b)$, so
$$g''(x/c-d/c)+cr(x)g'(x/c-d/c)+c^2s(x)g(x/c-d/c)=0,$$
and this equation is hypergeometric as well.
\end{proof}


\newpage

\begin{thebibliography}{99}

\bibitem{AnosovBolibruch}
D. V. Anosov and A. A. Bolibruch.  {\it The Riemann-Hilbert problem.}  Aspects of Mathematics E22.  Friedr. Vieweg \& Sohn, 1994.

\bibitem{Beauville}
A. Beauville.  {\it Monodromie des syst\`emes diff\'erentiels lin\'eaires \`a p\^oles simples sur la sph\`ere de Riemann (d'apr\`es A. Bolibruch).}  S\'eminaire Bourbaki, vol. 1992/93.  Ast\'erisque no. 216 (1993), exp. no. 765, 4, 103--119.

\bibitem{Belyi1}
G. V. Bely\u\i.  {\it Galois extensions of a maximal cyclotomic field.}  Izv. Akad. Nauk SSSR Ser. Mat. 43 (1979) 267--276, 479.

\bibitem{Belyi2}
G. V. Bely\u\i.  {\it A new proof of the three-point theorem.}  Mat. Sb. 193 (2002) 21--24; translation in Sb. Math. 193 (2002) 329--332.

\bibitem{BetremaZvonkin}
J. B\'etr\'ema and A. Zvonkin.  {\it Plane trees and Shabat polynomials.}  Proceedings of the 5th Conference on Formal Power Series and Algebraic Combinatorics (Florence, 1993).  Discrete Math. 153 (1996) 47--58.

\bibitem{GantzSteer}
C. Gantz and B. Steer.  {\it Gauge fixing for logarithmic connections over curves and the Riemann-Hilbert problem.}  J. London Math. Soc. (2) 59 (1999) 479--490.

\bibitem{Grothendieck}
A. Grothendieck. {\it Esquisse d'un programme.}  London Mathematical Society Lecture Note Series 242, Geometric Galois actions, 1, 5--48.  English translation on pp. 243--283.  Cambridge University Press, 1997. 

\bibitem{LandoZvonkin}
S. K. Lando and A. K. Zvonkin.  {\it Graphs on surfaces and their applications.}  With an appendix by Don B. Zagier.  Encyclopaedia of Mathematical Sciences 141.  Low-Dimensional Topology II.  Springer-Verlag, 2004. 

\bibitem{Maskit}
B. Maskit.  {\it Kleinian groups.}  Grundlehren der Mathematischen Wissenschaften 287.  Springer-Verlag, 1988.

\bibitem{Oesterle}
J. Oesterl\'e.  {\it Dessins d'enfants.}  S\'eminaire Bourbaki, vol. 2001/2002.  Ast\'erisque no. 290 (2003), exp. no. 907, ix, 285--305.

\bibitem{Remmert}
R. Remmert.  {\it Funktionentheorie I.}  Grundwissen Mathematik 5.  Springer-Verlag, 1984.

\bibitem{Schneps}
L. Schneps (editor).  {\it The Grothendieck theory of dessins d'enfants.}  Papers from the Conference on Dessins d'Enfants held in Luminy, April 19--24, 1993.  London Mathematical Society Lecture Note Series 200.  Cambridge University Press, 1994.

\bibitem{vanderPutSinger}
M. van der Put and M. F. Singer.  {\it Galois theory of linear differential equations.}  Grundlehren der Mathematischen Wissenschaften 328.  Springer-Verlag, 2003.

\end{thebibliography}
\end{document}